\documentclass[11pt,reqno]{amsart}

 \usepackage{amsmath}
 \usepackage{amssymb}
\usepackage{bm}
 \usepackage{mathtools}
\usepackage{color}
 
\usepackage{graphicx}

\newcommand{\mysection}[1]{\section{#1}
      \setcounter{equation}{0}}

\newtheorem{theorem}{Theorem}[section]
\newtheorem{lemma}[theorem]{Lemma}

\newtheorem{corollary}[theorem]{Corollary} 

\theoremstyle{definition}
\newtheorem{assumption}{Assumption}[section]

\theoremstyle{remark}
\newtheorem{remark}[theorem]{Remark} 

\newtheorem{example}{Example}[section]

\newcommand{\loc}{\text{\rm loc}}

 \makeatletter 
\makeatletter 
 \def\dashint{\operatorname{\,\,\,\mathclap{\int} \kern-.23em\text{\bf--}\!\!}}
 \makeatother

\def\dashnorm{\,\,\text{--}\kern-.4em\|}
\def\ninf{\qopname\relax\@empty{inf\phantom{p}\!\!\!}}

 \makeatother

\newcommand{\WO}{\overset{\scriptscriptstyle0}%
{W}\,\!}
  
\newcommand\bR{\mathbb{R}}

\newcommand\bC{\mathbb{C}}

\newcommand\bS{\mathbb{S}}

\newcommand\cL{\mathcal{L}}

\begin{document}

\title[Parabolic equations in Sobolev spaces]
{On parabolic equations in Sobolev spaces
with lower-order coefficients from   Morrey spaces}

\author{N.V. Krylov}
 
\email{nkrylov@umn.edu}
\address{127 Vincent Hall, University of Minnesota,
 Minneapolis, MN, 55455}


\dedicatory{Dedicated to N.N. Ural'tseva on the occasion of her 90'th birthday}

\begin{abstract}
We consider parabolic equations with operators
$\cL=\partial_{t}+a^{ij}D_{ij}+b^{i}D_{i}-c$ with $a$ being almost in VMO,
$b$ in a Morrey class containing $ L_{d+2}$  and $c $ in a Morrey class containing $L_{(d+2)/2}$. We prove
the solvability in Sobolev spaces of $\cL u=f\in L_{p}$ in bounded $C^{1,1}$-cylinders. 
\end{abstract}

\maketitle

\mysection{Introduction}

Let $\bR^{d}$, $d\geq 2$, be a Euclidean space of
points $x=(x^{1},...,x^{d})$.
 Fix a $\delta\in(0,1]$ and let $\bS_{\delta}$ be the set
of $d\times d$ symmetric matrices whose eigenvalues
are in $[\delta,\delta^{-1}]$.
Assume that on $\bR^{d+1}=\{(t,x):t\in\bR,x\in\bR^{d}\}$ we are given a Borel real-valued $c$, an $\bR^{d}$-valued 
function  $b=(b^{i})$, and an $\bS_{\delta}$-valued $a
=(a^{ij})$. Define
$$
\cL u=\partial_{t}u+a^{ij}D_{ij}u+b^{i}D_{i}u,
$$
where the summation over the repeated indices
is enforced and we use the notation
$$
 D_{i}  =\frac{\partial}{\partial x^{i}} ,\quad Du= (D_{i}u),\quad
  D_{ij} = D_{i}D_{j}  ,\quad D^{2}u=(D_{ij}u),
\quad
 \partial_{t} =\frac{\partial}{\partial t} .
$$

We are going
to   investigate
the equation
\begin{equation}
                                         \label{6.15.2}
 \cL u-(c+\lambda)u=f
\end{equation}
 in $\bR^{d+1}$, $(T,S)\times \bR^{d}$, and
$(T,S)\times \Omega$ with sufficiently large
$\lambda$. We are interested in Sobolev space
solutions and $b,c$ with rather poor integrability properties expressed in their belonging
to appropriate Morrey spaces.
 
In the last decades, much attention has been paid to regularity and solvability issues regarding boundary value problems for elliptic and parabolic operators with discontinuous coefficients. All this research was made possible thanks to the groundbreaking paper \cite{CFL_93}, where the Calder\'on--Zygmund theory of second-order elliptic operators with continuous principal coefficients has been extended to the discontinuous situation of operators with principal coefficients having vanishing mean oscillation
(VMO). Regularity and strong solvability theories have been developed for the Dirichlet, Neumann, and Robin boundary value problems both for elliptic and parabolic operators (see \cite{BC_93}, \cite{DK_11}, \cite{Kr_08},  \cite{MPS_00}, 
 \cite{PS_21} for instance). 
In \cite{ANPS_23} parabolic equations are considered
with Wentsell's boundary condition. Later
a different from \cite{CFL_93} method
was discovered
of treating equations with VMO coefficients,
which turned out to be applicable even to fully
nonlinear equations (see \cite{Kr_08}, \cite{Kr_18}).

In these works a substantial progress was achieved
in what concerns the regularity of the main coefficients. However, the assumptions on $b$ and $c$ remain, basically,
the same as in the classical monograph
\cite{LSU_67}: $b\in L_{d+2+}, c\in L_{(d+2)/2+}$. In the elliptic case, in the Sobolev space
theory, it was possible
to reduce the summability requirements on $b$ and $c$ by considering Morrey spaces
and using the results from \cite{Ad_86} or \cite{CF_90} (see, for instance, \cite{CLV_96},
\cite{Kr_23.1}, \cite{MT_13}, \cite{TTV_95}).

Recently in \cite{Kr_23} and \cite{8}
the author succeeded in obtaining   parabolic
analogues of some results in \cite{Ad_86} and \cite{CF_90} and this made it possible to
write this article on solvability in Sobolev spaces
of parabolic equations with singular coefficients.

The vast literature on equations solvable in Morrey
(mixed-norms) spaces is beyond the scope of this article.

We finish the introduction with some notation.
Define $B_{\rho}(x)=\{y \in\bR^{d} :|x-y|<\rho\}$, 
$$
C_{\rho}(t,x)=\big\{(s,y)\in \times\bR^{d +1 } :|x-y|<\rho,
t\leq s< t+\rho^{2}\big\},\quad C_{\rho}=C_{\rho}(0)
$$
and let $\bC_{\rho}$ be the collection of
$C_{\rho}(z)$, $z\in\bR^{d+1}$.
For measurable $\Gamma\subset \bR^{d+1}$
set $|\Gamma|$ to be its Lebesgue measure
and when it makes sense set
$$
 \dashint_{\Gamma}f\,dz=\frac{1}{|\Gamma|}
\int_{\Gamma}f\,dz.
$$
We fix $S,T$ such that
$-\infty\leq T<S\leq\infty$,   if $\Omega $
is a domain in $\bR^{d}$ we set $\Omega_{T,S}=(T,S)\times \Omega$.
By $W^{1,2}_{p}(\Omega_{T,S})$ we mean the collection of
functions $u$ such that $u,\partial_{t}u,Du,D^{2}u 
\in L_{p}(\Omega_{T,S})$. This space is provided
with a naturel norm. Its subset of functions
such that $u=0$ (in the well-known sense) on
$[[T,S]\times\partial\Omega]\cup[\{S\}\times\bar\Omega]$ is denoted by $\WO^{1,2}(\Omega_{T,S})$. In this notation we drop $S$
if $S=\infty$, in particular, $\bR^{d}_{0}=
(0,\infty)\times\bR^{d}$. Finally,  $L_{p}=L_{p}(\bR^{d+1})$.

\mysection{Main results}

Fix some $\rho_{0},\rho_{a}  \in(0,\infty)$
and let $1<p<q<\infty$.
Finite parameters $\theta$, $\check b$  and $\check c$  below will be specified later.
\begin{assumption}
                        \label{assumption 12.12.3}
 We have   
\begin{equation}
                                 \label{6.3.1}
 a^{\sharp}_{x,\rho_{a}}:=\sup_{\substack{\rho\leq\rho_{0}\\C\in\bC_{\rho}}}\dashint_{C}|a (t,x)-a_{C}(t)|
\,dx dt \leq  \theta,
\end{equation}
where
$$
a_{C}(t)=\dashint_{C}a (t,x)\,dxds\quad (\text{note}\quad t\quad \text{and}\quad ds).
$$
\end{assumption}

\begin{assumption}
                      \label{assumption 10.23.1}
 (i) If $p>d+2$, then
\begin{equation}
                           \label{10.23.1}
\sup_{C\in \bC_{\rho_{0}}} 
\| b \|_{L_{p}(C)}=:\bar b<\infty .
\end{equation}

 (ii) If $p=d+2$, then
\begin{equation}
                           \label{10.23.3}
\sup_{C\in \bC_{\rho_{0}}} 
\| b \|_{L_{q}(C)}=:\bar b<\infty .
\end{equation}

(iii) If $p<d+2$, then  
\begin{equation}
                           \label{3.14.2}
 \sup_{r\leq\rho_{0}}r
\sup_{C\in \bC_{r}} 
\dashnorm b \|_{L_{q}(C)}\leq   \check b .
\end{equation}
\end{assumption}

\begin{assumption}
                      \label{assumption 10.23.2}

(i) If $p>(d+2)/2$, then
\begin{equation}
                           \label{10.23.40}
\sup_{C\in \bC_{\rho_{0}}} 
\| c \|_{L_{p}(C)}=:\bar c<\infty .
\end{equation}

(ii) If $p=(d+2)/2$, then
\begin{equation}
                           \label{10.23.04}
\sup_{C\in \bC_{\rho_{0}}} 
\| c \|_{L_{q}(C)}=:\bar c<\infty .
\end{equation}

(iii) If $p<(d+2)/2$, then

\begin{equation}
                           \label{3.14.3}
 \sup_{r\leq\rho_{0}}r^{2}
\sup_{C\in \bC_{r}} 
\dashnorm c \|_{L_{q}(C)}\leq   \check c .
\end{equation}
\end{assumption}

\begin{remark}
                      \label{renark 10.21.1}
If $p<d+2<q$, in \eqref{3.14.2} one can take $q=d+2$.
Also if $p<(d+2)/2<q$, in \eqref{3.14.3} and 
in \eqref{3.14.2} we can take
$q=(d+2)/2$. This follows from H\"older's
inequality.
 
Also observe that, if $b,c$ are bounded,
Assumptions \ref{assumption 10.23.1} and \ref{assumption 10.23.2} are satisfied with any $\check b,\check c>0$ on account of taking $\rho_{0}$
small enough.
\end{remark}

\begin{remark}
                      \label{renark 10.25.1}
If 
$q\geq d+2$ the left-hand side of \eqref{3.14.2}
is
$$
N(d)\rho_{0}^{1-(d+2)/q}\sup_{C\in\bC_{\rho_{0}}}\|b\|_{L_{q}(C)},
$$
which can be made as small as we like on account of
choosing small $\rho_{0}$ if $b\in L_{q}$.
This also shows a significant improvement
over the conditions usually imposed in the Sobolev
space theory of parabolic equations
(see, for instance, Remark 3.3 in \cite{ANPS_23}).
\end{remark}

\begin{remark}
                           \label{remark 11.1.2}
There are examples when \eqref{3.14.2}
is satisfied and $b\not\in L_{q+\varepsilon,\loc}$,
no matter how small $\varepsilon>0$ is.
Similar statement holds concerning \eqref{3.14.3}.

\end{remark}

Here is our   main result.  
In the following theorem the constants $\check \lambda_{0}$
and $N$ depend on $\bar b$ only if $p\geq d+2$
and depend on $\bar c$ only if $p\geq (d+2)/2$.
Also the existence of $\check b$ ($\check c$)
is only needed if $p<d+1$ (respectively, $p<(d+2)/2$), but in order not to make the statement too heavy,
we list these  parameters in all cases anyway.
 
\begin{theorem}
                       \label{theorem 10.27.1}
 Let  $\Omega$
be a bounded $C^{1,1}$ domain in $\bR^{d}$
or $\Omega=\bR^{d}$.
Then under the above assumptions there exist 
$$
(\check b,\check c) =(\check b,\check c)(d,\delta, p, 
q, \rho_{a},\Omega  )\in(0,1]^{2},
$$
$$
\check \lambda_{0}=\check \lambda_{0}
(d,\delta, p, 
q, \rho_{a}, \rho_{0}, \bar b,\bar c ,\Omega)>0
$$
 such that,
 for  any
$\lambda\geq   \check \lambda_{0} $  and $f\in 
L_{p}(\Omega_{T,S})$ there exists
a unique $\WO^{1,2}_{p   }(\Omega_{T,S})$-solution $u$
of $\cL u-(c+\lambda) u=f$. Furthermore, there exists
a constant $N$ depending only on  $d$, $\delta$, $p$, $q$, $\rho_{a}$,  $\rho_{0}$,
  $\bar b$, $\bar c$,    $\Omega$,
such that 
\begin{equation}
                                \label{10.28.01}
\|\partial_{t} u,D^{2}u, \sqrt\lambda Du, \lambda u\|_{L_{p}(\Omega_{T,S})}
\leq   N   \|f\|_{L_{p}(\Omega_{T,S})}.    
\end{equation}
\end{theorem}

\begin{example}
                          \label{example 5.12.1}
It turns out that $\check b$ cannot be very large
in order for Theorem \ref{theorem 10.27.1} to hold.
To show this take $T=-1$, $S=0$, $\Omega=\bR^{d}$,
$d\geq2$,  
and for $\lambda>0$ set
$$
u(t,x)=e^{\lambda t+ |x|^{2}/(4t)}.
$$
This function satisfies
$$
\partial_{t}u+\Delta u-\frac{d}{|x|^{2}}x^{i}D_{i}u
 -\lambda u=0
$$
in $\Omega_{T,S}$. Also, as is easy to check, for $p <(d+2)/2$, we have
$ u,Du, D^{2}u,\partial_{t}u\in L_{p }(\Omega_{T,S})$ and $u\in\WO^{1,2}_{p}(\Omega_{T,S})$. 
 In addition, for $b=-xd/|x|^{2}$ the 
left-hand side of \eqref{3.14.2} is finite (and independent of
$\rho_{0}$) for $p<q <d$ (there are such $q$ since $d\geq2$).
However, estimate \eqref{10.28.01} fails,
which shows that the constant factor $d$ in $b$
is not sufficiently small.

\end{example}

For nonzero boundary value problem we have the
following result which is obtained from 
Theorem \ref{theorem 10.27.1} by well-known means
(cf. Remark \ref{remark 11.2.2}).
\begin{theorem}
                   \label{theorem 10.26.1}
Suppose $S-T<\infty$ and $\Omega$ is a bounded
$C^{1,1}$ domain in $\bR^{d}$. Then there exist 
$$
(\check b,\check c) =(\check b,\check c)(d,\delta, p, 
q, \rho_{a}, \Omega  )\in(0,1]^{2},
$$
such that for any $g\in W^{1,2}_{p}(\Omega_{T,S})$
there exists a unique $W^{1,2}_{p   }(\Omega_{T,S})$-solution $u$
of $\cL u-c u=f$ such that $u-g\in \WO^{1,2}_{p   }(\Omega_{T,S})$. Furthermore, 
there exists
a constant $N$ depending only on  $d$, $\delta$, $p$, $q$, $\rho_{a}$,  $\rho_{0}$,
  $\bar b$, $\bar c$,    $\Omega$, $S-T$,
such that 
\begin{equation}
                                \label{10.28.010}
\| u\|_{W^{1,2}_{p}(\Omega_{T,S})}
\leq   N   \|f\|_{L_{p}(\Omega_{T,S})}
+N\| g\|_{W^{1,2}_{p}(\Omega_{T,S})}.    
\end{equation}

\end{theorem}

Theorem \ref{theorem 10.27.1} is applicable
in the case of $a,b,c,f$   independent of $t$,
and owing to uniqueness yields solutions that
are also independent of $t$. This leads to
the ``elliptic'' counterpart
of our results, which however are much weaker than
in \cite{Kr_23.1} where in the assumptions we use $d$
in place of $d+2$.
 
\begin{remark}
                           \label{remark 11.1.1}
The following is slightly imprecise.
If $p<(d+2)/2$, the solutions in Theorem
\ref{theorem 10.27.1} are only of class
$W^{1,2}_{p}$ and thus not necessarily bounded.
However, if $f\in L_{(d+2)/2+\varepsilon}$,
$\varepsilon>0$, then $f$ is
in the Morrey class $E_{p,\beta}$ with $\beta<2$
and, according to \cite{5}, the solution is bounded.
This is the same situation as if there were no $b,c$.
Hence,
 our operator has Green's functions summable
to $(d+2)/d-\varepsilon$ power, despite very
low integrability of $b,c$.
\end{remark}

\mysection{Preliminary estimates}

                       \label{section 10.8.1}
 
For domains $Q\subset\bR^{d +1 }$, $p\in[1,\infty)$, and $\beta>0$, introduce
Morrey's space $E_{p,\beta}(Q)$
as the set of $g \in L_{p,\loc}(Q)  $ such that
\begin{equation}
                             \label{8.11.2}
\|g\|_{E_{p,\beta}(Q)}:=
\sup_{\rho<\infty,(t,x)\in Q}\rho^{\beta}
\dashnorm g I_{Q}\|_{L_{p}(C_{ \rho}(t,x) )} <\infty ,
\end{equation}
where
$$
\dashnorm g\|_{L_{p}(\Gamma)}=\Big(
\dashint_{\Gamma}|g|^{p}\,dz\Big)^{1/p}.
$$
We abbreviate $E_{p,\beta}=E_{p,\beta}(\bR^{d+1})$.
Note that, if $\beta p>d+2$, the space 
$E_{p,\beta}(Q)$ consists of only one function $g=0$.
Also observe that if $\beta\leq (d+2)/p$ and $Q=C_{R}$
one can restrict $\rho$ in \eqref{8.11.2}
to $\rho\leq R$.

For $s,r>0,\beta\in \bR$, and appropriate $f(t,x)$'s
on $\bR^{d+1}$ define
$$
p_{\beta}(s,r)=\frac{1}{s^{(d+2-\beta)/2}}e^{-r^{2}/s}I_{s>0},
$$
$$
P_{\beta}f(t,x)=\int_{\bR^{d }_{0} }p_{\beta}(s,|y|)f(t+s,x+y)\,dyds
$$
$$
=\int_{t}^{\infty}\int_{\bR^{d} }p_{\beta}(s-t,|y-x|)
f(s,y)\,dsdy.
$$
Observe that, if $f$ is independent of $t$,
then
$$
P_{\beta}f(t,x)=
P_{\beta}f(x)=N(\beta)\int_{\bR^{d}}
\frac{1}{|y|^{d-\beta}}f(x+y)\,dy= NI_{\beta}f(x),
$$
where $I_{\beta}$ is the Riesz potential.
In our investigation the most important values
of $\beta$ are 1 and 2.

 Here is a generalization of Theorem 4.3 of \cite{Kr_23} and which presents a parabolic
analog of a particular case of Theorem 7.3 of \cite{Ad_86} or, if $\beta=1$, of the Theorem
in \cite{CF_90}. Its proof is   different from
the arguments in \cite{Kr_23} and
\cite{Ad_86} but is based on the ideas from these
papers.
\begin{theorem}
                    \label{theorem 1.18.1}
Let $  q> p>1$, $\beta>0$,  $b\in E_{q,\beta}$.  Then for any $f\geq0$
we have
\begin{equation}
                         \label{5.30.1}
I:= \int_{\bR^{d+1}}|b|^{p}(P_{\beta}f)^{p}\,dz\leq
N\|b\|_{E_{q,\beta}}^{p}\|f\|_{L_{p}},
\end{equation}
where $N$ depends only on $d,p,q,\beta$.
In particular, 
 for any $u\in C^{\infty}_{0}$
\begin{equation}
                       \label{1.18.3}
 \int_{\bR^{d+1} }|b |^{p}|Du |^{p}\,dz
\leq N\|b\|_{E_{q,1}}^{p}K,
\quad
\int_{\bR^{d+1} }|b |^{p}| u |^{p}\,dz
\leq N\|b\|_{E_{q,2}}^{p}K,
\end{equation}
where
$
K=\|D^{2}u,\partial_{t}u\|^{p}_{L_{p}}
$
and $N$ depends only on $d,p,q$.

\end{theorem}
\begin{remark}
                     \label{remark 11.2.1}
The first estimate in \eqref{1.18.3} follows from \eqref{5.30.1} with $\beta=1$   and the fact that
for $f=\partial_{t}u+\Delta u$ we have
$$
Du(t,x)=c\int_{\bR^{d +1 }_{0}}
\frac{y}{s^{(d+2)/2}}e^{-|y|^{2}/(4s)}
f(t+s,x+y)\,dyds,
$$
where $c$ is a constant  and
$(|y|/s^{1/2})e^{-|y|^{2}/(4s)}\leq
Ne^{-|y|^{2}/(8s)}$. The second estimate
follows when $\beta=2$ since
$$
u(t,x)=c\int_{\bR^{d +1 }_{0}}
\frac{1}{s^{d/2}}e^{-|y|^{2}/(4s)}
f(t+s,x+y)\,dyds.
$$

\end{remark}

$$
M_{\gamma}f(t,x)=\sup_{\rho>0}\rho^{\gamma}
\dashint_{C_{\rho}(t,x)}|f(z)|\,dz,\quad
\gamma>0,
$$
$$
M f =M_{0}f .
$$

The next results might be a usefull complement
to Theorem \ref{theorem 1.18.1}. They are similar
to the corresponding results from \cite{Kr_23.1},
albeit with a restricted range of $p$ but with $q=p$.
We suppose that $p>2$ and do not know   if
the results still hold  for $p=2$. Recall that we write $b\in A_{1}$ if there is a constant $[b]_{A_{1}}$ such that $Mb \leq [b]_{A_{1}}
|b|$.
\begin{theorem}
                    \label{theorem 1.18.10}
Let $d+2\geq p >2$,   $b\in E_{p,1}, |b|^{p}\in A_{1}$.
Then for any $u\in C^{\infty}_{0}$
\begin{equation}
                       \label{1.18.03}
I:=\int_{\bR^{d+1} }|b |^{p}|Du |^{p}\,dz
\leq N\|b\|_{E_{p,1}}^{p}K,
\end{equation}
where
$
K=\|D^{2}u,\partial_{t}u\|^{p}_{L_{p}}
$
and $N$ depends only on $d,p$, and $[|b|^{p}]_{A_{1}}$.

\end{theorem}

\begin{theorem}
                    \label{theorem 11.5.1}
Let $d+2\geq p >2$,   $c\in E_{p,2}, |c|^{p}\in A_{1}$.
Then for any $u\in C^{\infty}_{0}$
\begin{equation}
                       \label{1.18.30}
I:=\int_{\bR^{d+1} }|c |^{p}| u |^{p}\,dz
\leq N\|c\|_{E_{p,2}}^{p}K,
\end{equation}
where
$
K=\|D^{2}u,\partial_{t}u\|^{p}_{L_{p}}
$
and $N$ depends only on $d,p$, and $[|c|^{p}]_{A_{1}}$.

\end{theorem}

\begin{remark}
                      \label{remark 11.5.1}
Usual arguments which one can find in
\cite{CF_90} or \cite{Kr_23.1} allow  one to
drop the conditions $|b|^{p}\in A_{1}$ and
$|c|^{p}\in A_{1}$ on account of requiring
$b\in E_{q,1}$ and $c\in E_{q,2}$ with
$d+2\geq q>p$. Of course then 
$\|b\|_{E_{p,1}}$ and $\|c\|_{E_{p,2}}$
shoud be replaced with $\|b\|_{E_{q,1}}$ and 
$\|c\|_{E_{q,2}}$, respectively.
\end{remark}

We prove these theorems after the proof of Theorem \ref{theorem 1.18.1} is given. For the latter
we need some preparations.

The following is, actually, a  corollary of the proof of part of Lemma 2.2 of \cite{Kr_23}.
We provided it with proof for completeness.

\begin{lemma}
                            \label{lemma 19.30.1}
For any $\alpha<\beta ,\beta\in(0,d+2]$, there exists a   constants $N$
($<\infty$)  such that for any $f\geq0$
  on $C_{1}$ we have
\begin{equation}
                             \label{1.17.2}
P_{\alpha}(I_{C^{c}_{2}}f)  
\leq N M_{\beta}
f(0)  .
\end{equation}
\end{lemma}

Proof. Clearly, we may assume that $f$ is bounded
with compact support.  Set $Q_{1}=\{(s,y):|y|\geq\sqrt s\}$,
$Q_{2}=\{(s,y):|y|\leq \sqrt s\}$.
Dealing with $P_{\alpha}(fI_{Q_{1}\cap C^{c}_{2}})$
we observe that $p_{\alpha}(s,r)\leq Nr^{-(d+2-\alpha)}$ since $\alpha\leq d+2$. Therefore,
$$
P_{\alpha}(fI_{Q_{1}\cap C^{c}_{2}})(0)
\leq N\int_{2}^{\infty}\frac{1}{r^{d+2-\alpha}}
\int_{0}^{r^{2}}\Big(\int_{|y|=r}f(s,y)\,d\sigma_{r}\Big)\,dsdr,
$$
where $d\sigma_{r}$ is the element of the surface area on $|y|=r$. By using that $\alpha<d+2$ and integrating by parts we get
$$
P_{\alpha}(fI_{Q_{1}\cap C^{c}_{2}})(0)\leq N
\int_{2}^{\infty}\frac{1}{r^{d+3-\alpha}}
\int_{2}^{r}\Big(\int_{0}^{\rho^{2}}
\Big(\int_{|y|=\rho}f(s,y)\,d\sigma_{\rho}\Big)\,ds\Big)\,d\rho dr
$$
$$
\leq N
\int_{2}^{\infty}\frac{1}{r^{d+3-\alpha}}
\int_{0}^{r}\Big(\int_{0}^{r^{2}}
\Big(\int_{|y|=\rho}f(s,y)\,d\sigma_{\rho}\Big)\,ds\Big)\,d\rho dr
$$
$$
=N\int_{2}^{\infty}\frac{1}{r^{d+3-\alpha}}I(r)\,dr ,
$$
where  
$$
I(r)=\int_{C_{r}}f(s,y)\,dyds.
$$
Observe that $r\geq 2$ and for any $(t,x)\in C_{1}(-1,0)$
$$
I(r)\leq 
\int_{C_{2r}(t,x)}f(s,y)\,dyds
$$
$$
\leq 
Nr^{ d+2-\beta }\Big((2r)^{\beta}\dashint_{C_{2r}(t,x)}f(s,y)\,dyds\Big)\leq Nr^{ d+2-\beta }M_{\beta}f(t,x).
$$
We also use that    $\alpha<\beta$. Then we see that 
\begin{equation}
                          \label{1.17.3}
P_{\alpha}(fI_{Q_{1}\cap C^{c}_{2}})(0)\leq N M_{\beta}f(t,x).
\end{equation}

Next, by integrating by parts
we  obtain  that
$$
P_{\alpha}(fI_{Q_{2}\cap C^{c}_{2}})(0)
\leq
\int_{4}^{\infty}\frac{1}{s^{(d+2-\alpha)/2}}
\int_{|y|\leq\sqrt s}f(s,y)\,dyds
$$
$$
\leq N\int_{4}^{\infty}\frac{1}{s^{(d+4-\alpha)/2}}I(\sqrt s)\,ds=
N\int_{2}^{\infty}\frac{1}{r^{d+3-\alpha}}I(r)\,dr.
$$
This along with \eqref{1.17.3} prove  that
$P_{\alpha}(I_{C^{c}_{2}}f)(0) 
\leq N M_{\beta}
f (t,x) $, which is equivalent to our statement.
The lemma is proved.
\begin{corollary}
                        \label{corollary 10.30.1}
There is a constant $N$ such that for any $f$,
for which $P_{\beta}|f|$ is finite on $C_{1}$,
on $C_{1}$ we have
\begin{equation}
                             \label{10.30.1}
|DP_{\beta}(I_{C^{c}_{2}}f) |+|\partial_{t}
P_{\beta}(I_{C^{c}_{2}}f)  |
\leq N M_{\beta}
f (0) .
\end{equation}
\end{corollary}

Indeed, one can concentrate on $f$ with compact
support 
when
$$
|D_{i}P_{\beta}f(t,x)|\leq\int_{\bR^{d+1} }
|D_{i}p_{\beta}(s,|y|)|\,|I_{C^{c}_{2}}f(t+s,x+y)|      \,dyds,
$$
where 
$$
|D_{i}p_{\beta}(s,|y|)|\leq N
\frac{1}{s^{(d+2-(\beta-1))/2}}e^{-r^{2}/(2s)}I_{s>0}.
$$
Similarly, for $(t,x)\in C_{1}$
$$
|\partial_{t}P_{\beta}f(t,x)|\leq\int_{\bR^{d+1} }
|\partial_{s}p_{\beta}(s,|y|)|\,|I_{C^{c}_{2}}f(t+s,x+y)|      \,dyds,
$$
where 
$$
|\partial_{s}p_{\beta}(s,|y|)|\leq N
\frac{1}{s^{(d+2-(\beta-2))/2}}e^{-r^{2}/(2s)}I_{s>0}.
$$

Corollary \ref{corollary 10.30.1} and the mean value theorem  yield
\begin{corollary}
                        \label{corollary 10.30.2}
There is a constant $N$ such that for any $f$,
for which $P_{\beta}|f|$ is finite on  $C_{1}$,
on $C_{1}$ we have
\begin{equation}
                             \label{10.30.10}
\int_{C_{1}} |
P_{\beta}(I_{C^{c}_{2}}f)(z )-P_{\beta}(I_{C^{c}_{2}}f)(0 )|\,dz 
\leq  N  M_{\beta}
f (0) .
\end{equation}
\end{corollary}

\begin{lemma}
                           \label{lemma 10.30.2}

There exists a constant $N$ such that for any $f$
\begin{equation}
                             \label{10.30.3}
\int_{C_{1}}P_{\beta}(I_{C _{2}}|f|) \,dz
\leq N  M_{\beta}f(0).
\end{equation}
\end{lemma}

Proof. We have
$$
\int_{C_{1}}P_{\beta}(I_{C _{2}}|f|)\,dz
=P_{\beta}g(0),
$$
where
$$
g(s,y)=\int_{C_{1}}(I_{C_{2}}|f|)(s+t,x+y)\,dxdt
\leq  I_{C_{3}}(s,y)M_{\beta}f(0)
$$
This proves \eqref{10.30.3} and the lemma.

The next result is similar to Proposition
3.3 of \cite{Ad_75} with a similar proof.
By the way, if $f$ is independent of $t$,
it yields Proposition
3.3 of \cite{Ad_75}.

\begin{theorem}
                           \label{theorem 10.30.1}
Let $\beta\in(0,d+2]$. Then there is a constant $N$ such that for any $f$,
for which $P_{\beta}|f|$ is
  locally summable,   we have
\begin{equation}
                                 \label{10.30.5}
P^{\sharp}_{\beta}f\leq N M_{\beta}f.
\end{equation}

\end{theorem}

Proof. It suffices to prove \eqref{10.30.5}
at the origin, that is, it suffice to prove that
for any $\rho>0$
\begin{equation}
                                 \label{10.30.6}
\dashint_{C_{\rho}}\dashint_{C_{\rho}}|
P_{\beta}f(z_{1})-P_{\beta}f(z_{2})|\,dz_{1}dz_{2}
\leq N M_{\beta}f(0).
\end{equation}
 By self-similarity, it suffices to concentrate
on $\rho=1$, in which case the left-hand side of
\eqref{10.30.6} is dominated by
$$
N\int_{C_{1}} |
P_{\beta}(I_{C^{c}_{2}}f)(z )-P_{\beta}(I_{C^{c}_{2}}f)(0)|\,dz 
+N\int_{C_{1}}P_{\beta}(I_{C _{2}}|f|) \,dz.
$$
After that the result follows immediately from
Corollary \ref{corollary 10.30.2} and Lemma \ref{lemma 10.30.2}. The theorem is proved.
  
By using the Fefferman-Stein theorem 
(Theorem 3.6 in \cite{GR_85} or Theorem 3.2.10 in \cite{Kr_08}) and the
Hardy-Littlewood maximal function theorem
we come to the following parabolic analogue
of the Muckenhoupt-Wheeden theorem
(see Theorem 3.6.1 in \cite{AH_96}).

\begin{corollary}
                       \label{corollary 10.30.4}
Let $r\in(1,\infty)$. Then there is a constant $N$
such that, for any $f\geq 0$,
$\|P_{\beta}f\|_{L_{r}}\leq N\|M_{\beta}f\|_{L_{r}}$.

\end{corollary}

{\bf Proof of Theorem \ref{theorem 1.18.1}}.
We may assume that $b$ and $f$ are bounded and have
compact support. This guarantees that $I<\infty$.
Then assume that $b\geq 0$, set
$u=P_{\beta}f$, and write   
\begin{equation}
                           \label{5.30.2}
I=\int_{\bR^{d+1}}\big(b^{p}u^{p-1}\big) P_{\beta}f
\,dz=\int_{\bR^{d+1}}P_{\beta}^{*}\big(b^{p}u^{p-1}\big)  f
\,dz\leq \|f\|_{L_{p}}\big\|
P_{\beta}^{*}\big(b^{p}u^{p-1}\big)\big\|_{L_{p'}},
\end{equation}
where $p'=p/(p-1)$ and $P^{*}_{\beta}$ is the conjugate operator to $P_{\beta}$, namely, for any $g\geq0$,   
\begin{equation}
                           \label{5.30.10}
(P^{*}_{\beta}g)
(s,x)=\big(P_{\beta}(g(-\cdot,-\cdot)\big)(-s,-x).
\end{equation}
 By introducing similarly $M_{\beta}^{*},M^{*}$
and using Corollary \ref{corollary 10.30.4} we conclude
$$
\|
P_{\beta}^{*}\big(b^{p}u^{p-1}\big)\big\|_{L_{p'}}
\leq N\|
M_{\beta}^{*}\big(b^{p}u^{p-1}\big)\big\|_{L_{p'}},
$$
where by H\"older's inequality
$$
M_{\beta}^{*}\big(b^{p}u^{p-1}\big)
=M_{\beta}^{*}\big(b(b^{p-1}u^{p-1})\big)\leq
\|b\|_{E_{q,\beta}} [M^{*}\big((bu)^{r}\big)]^{1/q'}
$$
with $r=q(p-1)/(q-1)$ and $q'=q/(q-1)$.
Since $p<q$, we have $p'>q'$ and by the Hardy-Littlewood theorem
$$
\|
M_{\beta}^{*}\big(b^{p}u^{p-1}\big)\big\|_{L_{p'}}
\leq N\|b\|_{E_{q,\beta}}\Big(\int_{\bR^{d+1}}
\big[(bu)^{r}\big]^{p'/q'}\,dz\Big)^{(p-1)/p}
$$
$$
=N\|b\|_{E_{q,\beta}}I^{(p-1)/p}.
$$
This, obviously, proves the theorem.
 
{\bf Proof of Theorem \ref{theorem 1.18.10}}. We follow the arguments in \cite{CF_90}.
 Define $V=P_{2}(|b|^{p})$,
so that $|b|^{p}=-N(\partial_{t}V+(1/4)\Delta V)$
with an appropriate $N$. It follows by integrating by
parts that 
$I=I_{1}+I_{2}$, where
$$
I_{1}=N\int_{\bR^{d+1}}  V|Du|^{p-2}D_{i}uD_{i}\partial_{t}u \,dz=-N\int_{\bR^{d+1}} D_{i}V|Du|^{p-2}D_{i}u \partial_{t}u\,dz
$$
$$
-N\int_{\bR^{d+1}}V(|Du|^{p-2}\Delta u+(p-2)|Du|^{p-4}D_{i}uD_{j}u
D_{ij}u\big)\partial_{t}u\,dz =:J_{1}+J_{2},
$$
$$
I_{2}=N\int_{\bR^{d+1}}D_{i}V
|Du|^{p-2}D_{j} u  D_{ij}u \,dz
$$ 
$$
\leq N\int_{\bR^{d+1}}P_{1}(|b|^{p})|Du|^{p-1}
|D^{2}u|\,dz.
$$

For $\alpha=2$ and $\beta=p$ ($>\alpha$) it follows from Lemma 2.2 of \cite{Kr_23} that 
$$
V\leq N\|b\|_{E_{p,1}}^{2}(M(|b|^{p}))^{1-2/p}
$$
 and, since
$M(|b|^{p})\leq\hat M(|b|^{p})\leq N_{0}|b|^{p}$, we have
$$
J_{2}\leq N\|b\|_{E_{p,1}}^{2}N_{0}^{1-2/p}\int_{\bR^{d+1}} |b|^{p-2} 
|Du|^{p-2}(|D^{2}u|^{2}+|\partial_{t}u|^{2})\,dxdt
$$
$$
\leq N\|b\|_{E_{p,1}}^{2}N_{0}^{1-2/p}
I^{(p-2)/p}\||D^{2}u|+|\partial_{t}u|\|^{2}_{L_{p}},
$$
where the last inequality is obtained by using
H\"older's inequality.

For $\alpha=1$ and $ \beta=p$ Lemma 2.2 of \cite{Kr_23} 
yields   
$$
P_{1}(|b|^{p})\leq N\|b\|_{E_{p,1}}
N_{0}^{1-1/p}|b|^{1-1/p}
$$ 
so that
$$
I_{2}\leq N\|b\|_{E_{p,1}}N_{0}^{1-1/p}
\int_{\bR^{d+1}}|b|^{p-1}|Du|^{p-1}
|D^{2}u| \,dz
$$
$$
\leq N\|b\|_{E_{p,1}}N_{0}^{1-1/p}
I^{(p-1)/p}\|D^{2}u \|_{L_{p}}.
$$
Similarly, $J_{1}$ is estimated and we conclude
that
$$
I\leq N \|b\|_{E_{p,1}}^{2}N_{0}^{1-2/p}
I^{(p-2)/p} K^{2/p}+N\|b\|_{E_{p,1}}N_{0}^{1-1/p}
I^{(p-1)/p} K^{1/p}.
$$
Here $I<\infty$, since $u$ has compact support, and
in terms of  $\mu:= K\|b\|^{p}_{E_{p,1}}/(N_{0}I)$
the above  inequality means that
$N_{0}^{-1}\leq N\mu^{2/p}+N\mu^{1/p}$,
so that either $\mu\leq 1$ and $N\mu^{1/p}
\geq N_{0}^{-1}$, $NK\|b\|^{p}_{E_{p,1}}
\geq IN_{0}^{1-p} $, or $\mu\geq 1$
and $I\leq N_{0}^{-1}K\|b\|^{p}_{E_{p,1}}$. Since $N_{0}\geq1$,
in both cases

\begin{equation}
                              \label{1.18.5}
I\leq NN_{0}^{p-1}\|b\|^{p}_{E_{p,1}}K.
\end{equation}
The theorem is proved.

{\bf Proof of Theorem \ref{theorem 11.5.1}}.
Introduce $V$ in the same way as in the proof of
Theorem \ref{theorem 1.18.10} by replacing
$|b|$ with $|c|$ and assume that $c\geq0$. Then we get that $I=I_{1}+I_{2}+I_{3}$,
where
$$
I_{1}=N\int_{\bR^{d+1}}V|u|^{p-2}u\partial_{t}u\,dz,
\quad
I_{2}=N\int_{\bR^{d+1}} V|u|^{p-2}u\Delta u\,dz,
$$
$$
I_{3}=N\int_{\bR^{d+1}} V|u|^{p-2}|D u|^{2}\,dz.
$$

By Lemma 2.2 of \cite{Kr_23} with $\alpha=2$,
$\beta=2p$ we have 
$$
V\leq N\|c\|_{E_{p,2}}\big(M(c^{p}\big)^{1-1/p}
\leq N\|c\|_{E_{p,2}}c^{p-1}.
$$
This shows that 
$$
|I_{1}|+|I_{2}|\leq N\|c\|_{E_{p,2}}I^{1-1/p}K^{1/p}.
$$

While estimating $I_{3}$ we use that $c^{p-1}|u|^{p-2}
|Du|^{2}=(c^{p-2}|u|^{p-2})(c|Du|^{2})$ and apply
H\"older's inequality. Then we get
$$
|I_{3}|\leq NI^{1-2/p}\Big(\int_{\bR^{d+1}}
c^{p/2}|Du|^{p}\,dz\Big)^{2/p}
\leq NI^{1-2/p}\|c\|_{E_{p,2}}K^{2/p},
$$
where the last inequality follows from
Theorem \ref{theorem 1.18.10} since
$c^{1/2}\in E_{p,1}$ and $\|c^{1/2}\|_{E_{p,1}}
\leq \|c \|_{E_{p,2}}$. After that it only
remains to repeat the end of the proof
of Theorem \ref{theorem 1.18.10}. The theorem
is proved.

The main fact following from  Theorem \ref{theorem 1.18.1} is its version
in which the condition on $b$ is ``local''.

\begin{theorem}
                      \label{theorem 10.21.1}
Under Assumption \ref{assumption 10.23.1}, if $p<d+2$,
then for any $u\in C_{0}^{\infty}$,
\begin{equation}
                            \label{10.21.1}
\int_{\bR^{d+1} }|b |^{p}|Du |^{p}\,dz
\leq N\check b^{p}K,
\end{equation}
where
$
K=\|D^{2}u,\partial_{t}u\|^{p}_{L_{p}}
+\rho_{0}^{-2p}\| u\|^{p}_{L_{p}}
$
and $N$ depends only on $d,p,q$.

\end{theorem}

Proof.  H\"older's inequality allows us to
assume that $q\leq d+2$ in \eqref{3.14.2}.
Then take $\zeta,\eta\in C^{\infty}_{0}(C_{\rho_{0}})$,
$\zeta,\eta\geq0$, such that $\eta= \rho_{0}^{-(d+2)/p}$ on the support of $\zeta$ and
\begin{equation}
                         \label{11.4.5}
\int_{C_{\rho_{0}}}\zeta^{ p}\,dz=1, \quad \zeta+
\eta+ \rho_{0}|D\eta|
+ \rho^{2}_{0}(|D^{2}\eta|+|
\partial_{t}\eta|)\leq N(d)\rho_{0}^{-(d+2)/p}.
\end{equation}
 We claim that for any $\rho>0$ and  
$C\in \bC_{\rho}$   we have
\begin{equation}
                              \label{11.4.4}
\Big(\dashint_{C} |b \zeta|^{q}\,dz\Big)^{1/q}\leq N\rho_{0}^{-(d+2)/p}\check b \rho^{-1}.
\end{equation}
Indeed, if $\rho\leq \rho_{0}$ it suffices
to observe that $\zeta\leq N\rho_{0}^{-(d+2)/p}$.
In case $\rho>\rho_{0}$, it suffices to use that
$$
\dashint_{C } |b \zeta|^{q}\,dz=
N\rho^{-(d+2)}\int_{C }|b \zeta |^{q}\,dz
\leq 
N\rho_{0}^{-q(d+2)/p}\rho ^{-(d+2)}\int_{C_{\rho_{0}}}|b    |^{q}\,dz 
$$
$$
=N\rho_{0}^{-q(d+2)/p}\rho_{0}^{d+2}\rho ^{-(d+2)}
\dashint_{B_{\rho_{0}}}|b    |^{q}\,dz
\leq N\rho_{0}^{-q(d+2)/p}(\rho_{0}/\rho)^{d+2} \check b^{q}
\rho_{0}^{-q}
$$
$$
=N\rho_{0}^{-q(d+2)/p}(\rho_{0}/\rho)^{d+2-q}\rho^{-q}
\check b^{q}\leq N\rho_{0}^{-q(d+2)/p}\check b^{q}\rho^{-q},
$$
where at the last step we used that $q\leq d+2$.
Now, in light of \eqref{11.4.4} by   Theorem
  \ref{theorem 1.18.1}
 $$
\int_{\bR^{d+1}}|b \zeta|^{p}|Du|^{p}\,dz
\leq N\rho_{0}^{-(d+2) }\check b^{p}
\int_{\bR^{d+1}}(|D^{2} u|+|\partial_{t}u|)^{p}\,dz,\quad u\in C^{\infty}_{0}.
$$
Since $\eta$ is constant on the support of $\zeta$ and this inequality is also valid for $\eta u$ in place of $u$, we have
$$
\int_{\bR^{d+1}}|b \zeta|^{p}|Du|^{p}\,dz
\leq N \check b^{p}
\int_{\bR^{d+1}}(|D^{2} (\eta u)|+|\partial_{t}
(\eta u)|)^{p}\,dz,\quad u\in C^{\infty}_{0}.
$$

We plug in here $\zeta_{y}$ and
$\eta_{y}u$ in place of $\zeta$ and $\eta$, respectively, where $\zeta_{y}(z)=\zeta(z+y)$,
$\eta_{y}(z)=\eta(z+y)$. Then we get
$$
\int_{\bR^{d+1}}\zeta^{ p}_{y}|b |^{p}
|Du|^{p}\,dz\leq
N \check b^{p}
\int_{\bR^{d+1}}\eta^{p}_{y}(|D^{2}u| +
|\partial_{t}u|)^{p}\,dz
$$
$$
+N \check b^{p}
\int_{\bR^{d+1}}\big(|D\eta_{y}| |D u| 
+(|D^{2}\eta_{y}|+|\partial_{t}\eta_{y}||u|\big)^{p}\,dz.
$$
After integrating through with respect to $y$
and using that by H\"older's inequality 
and \eqref{11.4.5}
$$
\int_{\bR^{d+1}}\eta_{y}^{p}\,dy\leq
N ,\quad
\int_{\bR^{d+1}}|D\eta_{y}  |^{p}\,dy\leq
N\rho_{0}^{ -p},
$$
$$
\int_{\bR^{d+1}}(|D^{2}\eta_{y}  |+|\partial_{t}\eta_{y}|)^{p}\,dy\leq
N\rho_{0}^{ -2p}
$$
we come to \eqref{10.21.1} with $K$ also containing
the term
 $
 \rho_{0}^{-p}\|Du\|^{p}_{L_{p}}
 $, which is taken care of by interpolation.
This proves the theorem.

\begin{theorem}
                      \label{theorem 10.21.10}
Under Assumption \ref{assumption 10.23.2}, if $p<(d+2)/2$,
then for any $u\in C_{0}^{\infty}$,
\begin{equation}
                            \label{10.21.2}
\int_{\bR^{d+1} }|c |^{p}| u |^{p}\,dz
\leq N{\check c}^{p}K,
\end{equation}
where
$
K=\|D^{2}u,\partial_{t}u\|^{p}_{L_{p}}
+\rho_{0}^{-2p}\| u\|^{p}_{L_{p}}
$
and $N$ depends only on $d,p,q$.

\end{theorem}

Proof.  H\"older's inequality allows us to
assume that $q\leq (d+2)/2$ in \eqref{3.14.3}.
Then take the same $\zeta,\eta$ as in the previous proof
and observe that for $C\in \bC_{\rho}$
\begin{equation}
                              \label{11.4.40}
\Big(\dashint_{C} |c \zeta|^{q}\,dz\Big)^{1/q}\leq N\rho_{0}^{-(d+2)/p}\check c \rho^{-2}.
\end{equation}
This is proved in the same way as \eqref{11.4.4}.
After that applying Theorem \ref{theorem 1.18.1}
we finish the proof as in the case of Theorem
\ref{theorem 10.21.1}. The theorem is proved.

\begin{theorem}
                        \label{theorem 10.24.1}

Let Assumptions \ref{assumption 10.23.1}
and \ref{assumption 10.23.2} be satisfied,
let $\rho\in(0,\infty]$ and $u\in C^{\infty}_{0}$.
Then

(i) For $p< d+2$ there exists $N=N(d,p,q )$
such that
\begin{equation}
                              \label{10.23.6}
\|\,|b|\,|Du|\|_{L_{p}(C_{\rho})}\leq N 
 \check b
\|D^{2}u\|_{L_{p}(C_{\rho})}
+N \check b(\rho^{-2}+ \rho_{0}^{-2} )\| u\|_{L_{p}(C_{\rho})}.
\end{equation}

(ii) For $p< (d+2)/2$ there exists $N=N(d,p,q )$
such that
\begin{equation}
                              \label{10.24.1}
\|c  u \|_{L_{p}(C_{\rho})}\leq N 
 \check c 
\|D^{2}u\|_{L_{p}(C_{\rho})}
+N \check c(\rho^{-2}+ \rho_{0}^{-2 } )\|  u\|_{L_{p}(C_{\rho})}.
\end{equation}

\end{theorem}

Proof. 
(i) Estimate \eqref{10.23.6} for  
 $\rho=\infty$ is proved in Theorem \ref{theorem 10.21.1}.  If $\rho<\infty$, take any
extension operator $\Pi=\Pi_{\rho}$ which extends
smooth functions in $C_{\rho}$ to smooth
functions with compact support in $\bR^{d+1}$ (say, first extending them to $\bR\times B_{\rho}$ and then further) and is such that
$$
\|\Pi u\|_{L_{p}}\leq N\|u\|_{L_{p}(C_{\rho})},
\quad \|D\Pi u\|_{L_{p}}\leq N\|Du\|_{L_{p}(C_{\rho})}+N\rho^{-1}\|u\|_{L_{p}(C_{\rho})},
$$
$$
\|D^{2}\Pi u\|_{L_{p}}\leq N\|D^{2}u\|_{L_{p}(C_{\rho})}+N\rho^{-2}\|u\|_{L_{p}(C_{\rho})},
$$
 where $N=N(d,p)$. By the way,
the fact that $N$ can be chosen independent of $\rho$
is easily proved by rescaling. Then after applying
\eqref{10.23.6} with $\rho=\infty$ to $\Pi u$
we obtain \eqref{10.23.6} as is.

Assertion (ii) is proved similarly on the basis of 
Theorem \ref{theorem 10.21.10}. The theorem is proved.

\begin{remark}
                          \label{remark 3.3.1}
Theorem \ref{theorem 10.24.1} will still hold if
we replace cylinders $C_{\rho}$ with half-cylinders
(the base split in half). To see this it suffices to just
extend our functions across the flat part to the whole cylinder. Actually the boundary of ``half cylinders'' even need not to be flat, as long as it allows one to extend the functions $u$
across the border to $\hat u$ not much distorting the $L_{p}$ norms of $u,Du,D^{2}u,\partial_{t}u$. Therefore we can consider
$C_{\rho}(t,x)\cap \Omega_{T}$, where $\Omega$ is a bounded domain of class $C^{1,1}$,  
 $t\in (T,\infty)$, $x\in\partial \Omega$. Of course, in this situation
$\rho$ should be sufficiently small,
$\rho=\rho(d,p,\rho_{0},\Omega)$. However,
having it small enough, we can have
$$
\|\hat u\|_{L_{p}(C_{\rho}(t,x))}\leq N_{0}
\| u\|_{L_{p}(C_{\rho}(t,x)\cap\Omega_{T})},
$$
$$
\|D\hat u\|_{L_{p}(C_{\rho}(t,x))}\leq N_{0}
\|D u\|_{L_{p}(C_{\rho}(t,x)\cap\Omega_{T})}
+N_{1}\| u\|_{L_{p}(C_{\rho}(t,x)\cap\Omega_{T})}, 
$$
$$
\|D^{2}\hat u\|_{L_{p}(C_{\rho}(t,x))}\leq N_{0}
\|D^{2} u\|_{L_{p}(C_{\rho}(t,x)\cap\Omega_{T})}
+N_{1}\| u\|_{L_{p}(C_{\rho}(t,x)\cap\Omega_{T})}, 
$$
where $N_{0}=N(d,p)$, $N_{1}=N_{1}(d,p,\rho_{0},\Omega)$. This and partitions
of unity lead to (i) and (ii) in the following result
if $\Omega$ is a bounded domain of class $C^{1,1}$.
If $\Omega=\bR^{d}$ it suffices to take
$C_{1}(t,x),1$ in place of $C_{\rho},\rho$
in Theorem \ref{theorem 10.24.1} and then integrate
with respect to $(t,x)$ over $\Omega_{T}$.

\end{remark}

\begin{theorem}   
                         \label{theorem 3.3.1}
Suppose  that   Assumptions \ref{assumption 10.23.1}
and \ref{assumption 10.23.2} are satisfied,  $\Omega$ is a bounded domain
in $\bR^{d}$ of class $C^{1,1}$ or $\Omega=\bR^{d}$ and
let  $u\in W^{1,2}_{p}(\Omega_{T})$.

(i)
For $p<d+1 $ there exist  constants $N_{0}=N_{0}(p,q ,  d,\Omega )$ and $N_{1}=N_{1}(p,q ,  d, \rho_{0},\Omega )$ such that
$$
\|\,|b|\,|Du|\|_{L_{p}(\Omega_{T})}\leq N_{0}
 \check b
\|D^{2}u,\partial_{t}u\|_{L_{p}(\Omega_{T})}
+N_{1}\check b\|  u\|_{L_{p}(\Omega_{T})}.
$$

(ii) For $ p< (d+2)/2$ there exist constants $N_{0}=N_{0}(p,q ,  d ,\Omega)$ and $N_{1}=N_{1}(p,q ,  d,\rho_{0},\Omega )$ such that 
$$
\|cu\|_{L_{p}(\Omega_{T})}\leq N_{0}\check c
 \|D^{2}u,\partial_{t}u\|_{L_{p}(\Omega_{T})}+N_{1}\check c
\| u\|_{L_{p}(\Omega_{T})} .
$$

(iii) For $p\geq d+2$, for any $\varepsilon>0$,   there exists   
$N_{b}=N(\varepsilon, p,q ,d,\rho_{0},\bar b, \Omega)$ 
such that
\begin{equation}
                                \label{10.24.5}
\||b|\,|Du| \|_{L_{p}(\Omega_{T})}\leq 
\varepsilon\|D^{2}u,\partial_{t}u\| _{L_{p}(\Omega_{T})}+N_{b}
\| u\| _{L_{p}(\Omega_{T})} .
\end{equation}

(iv) For $p\geq (d+2)/2$, for any $\varepsilon>0$,   there exists   $N_{c}=N(\varepsilon, p,q ,d,\rho_{0},\bar c, \Omega)$
such that
\begin{equation}
                                \label{10.24.50}
\|cu\|_{L_{p}(\Omega_{T})}
\leq  
\varepsilon\|D^{2}u,\partial_{t}u\| _{L_{p}(\Omega_{T})}+N_{c}
\| u\| _{L_{p}(\Omega_{T})} .
\end{equation} 
 
\end{theorem}

Assertions (iii) and (iv) are classical
consequences of interpolation and embedding theorems.

\begin{remark}
                            \label{remark 11.2.2}
 A simple consequence of Theorem
\ref{theorem 3.3.1} is that $\cL$ is a bounded operator from $W^{1,2}_{p}(\Omega_{T})$ into
$L_{p}(\Omega_{T})$. In particular, the problem
of solving $(\lambda-\cL)u=f\in L_{p}(\Omega_{T})$ in $W^{1,2}_{p}(\Omega_{T})$ with boundary condition
$u-g \in \WO^{1,2}_{p}(\Omega_{T})$, where
$g\in W^{1,2}_{p}(\Omega_{T})$, reduces to
solving  $(\lambda-\cL)w=h\in L_{p}(\Omega_{T})$ in $\WO^{1,2}_{p}(\Omega_{T})$ by using the substitution
$w=u-g$, $h=f-(\lambda-L)g$.
\end{remark}

 \mysection{Proof of Theorem \protect\ref{theorem 10.27.1}}

Here is the combination of Theorems 2 and 6 of \cite{DK_11}, specified
to the case of one second order equation,
in which $\Omega=\bR^{d}$ or 
$\Omega$ is a bounded domain
of class $C^{1,1}$. 
Set
$$
\cL_{0}=\partial_{t}+a^{ij}D_{ij}.
$$

\begin{theorem}
                   \label{theorem 10.26.3}
There exist $\theta=\theta
(d,\delta,p)>0$ and $\lambda_{0}=
\lambda_{0}(d,\delta,p,\rho_{a},\Omega)>0$ such that,
if Assumption \ref{assumption 12.12.3} is satisfied with the above $\theta$, then for any $f\in L_{p}(\Omega_{T})$
and $\lambda\geq \lambda_{0}$
there exists a unique $u\in \WO^{1,2}_{p}(\Omega_{T})$
such that $\cL_{0}u-\lambda u=f$ in $\Omega_{T}$.
Furthermore,
$$
\|\partial_{t}u,D^{2}u,\sqrt \lambda Du,\lambda u
\|_{L_{p}(\Omega_{T})}\leq N\|f\|_{L_{p}(\Omega_{T})},
$$
where $N$ depends only on $d,\delta,p,\rho_{a},
\Omega$.
\end{theorem} 

{\bf Proof of Theorem  \ref{theorem 10.27.1}}.
One can assume that $S=\infty$. Indeed, if $S<\infty$
one can set $f(t,x)=0$ for $t\geq S$, find the solution  in $\Omega_{T}$, which will also be a solution in $\Omega_{S}$ and by uniqueness will
vanish for $t\geq S$ and, hence, belong to $\WO^{1,2}
(\Omega_{T,S})$. Dealing with $\WO^{1,2}
(\Omega_{T })$ and having in mind the method of continuity we convince ourselves that it suffices
to prove \eqref{10.28.01} as an a priori
estimate (with $S=\infty$). If $u\in \WO^{1,2}
(\Omega_{T })$ and $\cL u-(c+\lambda) u=f$,
then by Theorem \ref{theorem 10.26.3} for 
$\lambda\geq\lambda_{0}(d,\delta,p,\rho_{a},\Omega)$
\begin{equation}
                            \label{10.31.1}
\|\partial_{t}u,D^{2}u,\sqrt \lambda Du,\lambda u
\|_{L_{p}(\Omega_{T})}\leq \bar N \|f\|_{L_{p}(\Omega_{T})}
+\bar N \|b^{i}D_{i}u,cu\|_{L_{p}(\Omega_{T})},
\end{equation}
where $\bar N =\bar N (d,\delta,p,\rho_{a},\Omega)$.
The remaining part of the proof splits into a few cases
out of which we are going to do only one: $p<(d+2)/2$,
since the remaining cases are dealt with quite
similarly. In this case by Theorem
\ref{theorem 3.3.1} there exists   constants
$N_{0}=N_{0}(p,q,d,\Omega)$ and $N_{1}=N_{1}(p,q ,  d,\rho_{0},\Omega )$ such that
$$
\|b^{i}D_{i}u,cu\|_{L_{p}(\Omega_{T})}
\leq N_{0}(\check b+\check c)
\|D^{2}u,\partial_{t}u\|_{L_{p}(\Omega_{T})}+
N_{1}(\check b+\check c)\|u\|_{L_{p}(\Omega_{T})}.
$$
We now choose $\check b+\check c$,
depending only on $d,\delta,p,q,\rho_{a},\Omega$ so that
$$
\bar N N_{0}(\check b+\check c)\leq 1/2.
$$
This will allow us to absorb the term with $\|D^{2}u,\partial_{t}u\|_{L_{p}(\Omega_{T})}$ coming
from the last term in \eqref{10.31.1} into the left-hand side of \eqref{10.31.1}. After that it only remains to increase $\lambda_{0}$ to
$\lambda_{1}(d,\delta,p,q,\rho_{a},\rho_{0},\Omega)$
to absorb $\bar NN_{1}(\check b+\check c)\|u\|_{L_{p}(\Omega_{T})}$ into the left-hand side of
\eqref{10.31.1}. The theorem is proved.

{\bf Acknowledgments}. The author is sincerely grateful
to D. Kinzebulatov for pointing out
paper \cite{Ad_86} and to A. Lerner  who
communicated the source of
Theorem \ref{theorem 10.30.1}, which is
Proposition 3.3 of \cite{Ad_75}.


\begin{thebibliography}{mm}

 \bibitem{Ad_75} D. Adams,   A note on Riesz potentials,
{\em Duke Math. J.\/}, Vol. 42 (1975), No. 4, 765-778, 
doi.org/10.1215/s0012-7094-75-04265-9

\bibitem{Ad_86} D. Adams, {\em Weighted
nonlinear potential theory\/}, Trans. Amer. Math.
Soc., Vol. 297 (1986), No. 1, 73--94.

\bibitem{AH_96} 
D.R. Adams and L.I. Hedberg,   ``Function spaces and potential theory'', Grundlehren der mathematischen Wissenschaften [Fundamental Principles of Mathematical Sciences], 314. Springer-Verlag, Berlin, 1996. xii+366 pp.

\bibitem{ANPS_23} D.E. Apushkinskaya,
A. I.  Nazarov, D. K. Palagachev, L. G. Softova,  {\em Nonstationary Venttsel problems with discontinuous data\/},  J. Differential Equations, Vol. 375 (2023), 538--566.

\bibitem{BC_93} M. Bramanti and M.C. Cerutti, {\em $W^{1,2}_{p}$-solvability for the Cauchy-Dirichlet problem for parabolic equations with VMO coefficients\/}, Commun. Partial Differ. Equ. Vol. 18 (1993),
No.  9--10,  1735--1763.

\bibitem{CLV_96} P. Cavaliere, M. Longobardi,
and A. Vitolo, {\em Imbedding estimates and elliptic
equations with discontinuous coefficients in unbounded domains\/}, Le Matematiche, Vol. 51 (1996), No. 1,
87--104.

\bibitem{CF_90} F. Chiarenza and M. Frasca,
A remark on a paper by C. Fefferman,
{\em Proc. Amer. Math. Soc.\/}, Vol. 108 (1990), No. 2,
407--409, doi.org/10.2307/2048289

\bibitem{CFL_93}
F. Chiarenza, M. Frasca, P. Longo, {\em $W^{2,p}$--solvability of the Dirichlet problem for nondivergence elliptic equations with VMO coefficients\/}, Trans. Am. Math. Soc., Vol. 336 (1993), No.  2,  841--853.

\bibitem{DK_11} Hongjie Dong and Doyoon Kim,
{\em On the $L_{p}$-solvability of higher order parabolic and elliptic system with BMO coefficients\/},
Arch. Rational Mech. Anal., Vol. 199 (2011), 889--941.



\bibitem{FHS_17} G. Di Fazio, D.I. Hakim,
and Y. Sawano, Elliptic equations with discontinuous coefficients in generalized
Morrey spaces, {\em Europ. J. Math.\/}, Vol, 3 (2017), 729--762, doi.org/10.1007/s40879-017-0168-y



\bibitem{GR_85} J. Garcia-Guevra and J.L. Rubio de Francia, Weighted norm inequalities and related
topics, North-Holland Mat. Stud., Vol. 116, 1985, ISBN: 0 444 87804 1


\bibitem{Kr_08} N.V. Krylov,
``Lectures on elliptic and parabolic equations
in Sobolev spaces'', Amer.
Math. Soc., Providence, RI, 2008, doi.org/10.1090/gsm/096/04

 \bibitem{Kr_18} N.V. Krylov,
``Sobolev and viscosity solutions for fully nonlinear  elliptic 
and parabolic equations'', Mathematical Surveys and Monographs,
233, Amer.
Math. Soc., Providence, RI, 2018.

\bibitem{Kr_23} N.V. Krylov,  {\em On parabolic Adams's, the  Chiarenza-Frasca
theorems,
and some other results related to parabolic Morrey spaces\/}, Mathematics in Engineering, 
Vol. 5 (2023), No. 2, Paper No. 038, 20 pp.

\bibitem{Kr_23.1} N.V. Krylov, {\em Elliptic  equations in Sobolev spaces with
Morrey drift and the zeroth-order coefficients\/}, 
Trans. Amer. Math. Soc., Vol. 376 (2023), No. 10,  
 7329--7351
http://arxiv.org/abs/2204.13255

\bibitem{5} N.V. Krylov,  {\em On parabolic equations in Morrey spaces with
  VMO $a$ and Morrey $b,c$\/},  http://arxiv.org/abs/2304.03736

\bibitem{8} {\em A remark on a paper of
F. Chiarenza and M. Frasca\/},  
http://arxiv.org/abs/2310.12170

\bibitem{LSU_67}
O.A. Ladyzhenskaya, V.A. Solonnikov, N.N. Ural'tseva,
``Linear and quasi-linear parabolic equations'',
Nauka, Moscow, 1967 in Russian: English translation: American Math. Soc.,
Providence, 1968.

\bibitem{MPS_00}
A. Maugeri, D.K. Palagachev, and L.G. Softova, ``Elliptic and Parabolic Equations with Discontinuous Coefficients'', Mathematical Research, Vol.109, Wiley-VCH Verlag Berlin GmbH, Berlin, 2000.

\bibitem{MT_13} S. Monsurr\`o and M. Transirico,
{\em A priori bounds in $L^{p}$ and $W^{2,p}$
for solutions of elliptic equations\/}, Abstract
and Applied Analysis, Volume 2013, Article ID 650870,
7 pages, http://dx.doi.org/10.1155/2013/650870

\bibitem{PS_21} D.K. Palagachev and L.G. Softova,
{\em Generalized Morrey regularity of $2b$-parabolic systems\/}, Applied Math. Letters, Vol/ 112 (2021),106838.

\bibitem{TTV_95} M.Transirico, M.Troisi, and A.Vitolo,
{\em Spaces of Morrey type and elliptic equations in 
divergence form on unbounded domains\/},
 Bollettino della Unione Matematica Italiana,
(7) 9-B (1995), 153--174.
 
\end{thebibliography}
\end{document}